%% file: main.tex
\documentclass[11pt,reqno]{amsart}
\usepackage{amssymb,blkarray}
\usepackage[colorlinks]{hyperref}
\usepackage[all]{xy}

\title{Mukai flops and $\P$-twists}

\author[N.~Addington]{Nicolas Addington}
\address{Nicolas Addington \\
Department of Mathematics \\
University of Oregon \\
Eugene, OR 97403-1222 \\
United States}
\email{adding@uoregon.edu}

\author[W.~Donovan]{Will Donovan}
\address{Will Donovan \\
Kavli Institute for the Physics and Mathematics of the Universe (WPI) \\
The University of Tokyo Institutes for Advanced Study \\
5-1-5 Kashiwanoha \\
Kashiwa, Chiba, 277-8583 \\
Japan}
\email{will.donovan@ipmu.jp}

\author[C.~Meachan]{Ciaran Meachan}
\address{Ciaran Meachan \\ School of Mathematics \\
The University of Edinburgh \\
James Clerk Maxwell Building \\
Peter Guthrie Tait Road \\
Edinburgh EH9 3FD \\
United Kingdom}
\email{ciaran.meachan@ed.ac.uk}

\newlength{\baselineskipplus}
\newcommand\xymatrixresized{\xymatrix@R=\baselineskip@C=\baselineskipplus}

\usepackage{versions}
\includeversion{cut}
\excludeversion{ins}

\newcommand \A {\mathbb A}
\newcommand \C {\mathbb C}
\newcommand \E {\mathcal E}
\newcommand \F {\mathcal F}
\newcommand \G {\mathcal G}
\renewcommand \H {\mathcal H}
\renewcommand \L {\mathcal L}
\newcommand \M {\mathcal M}
\newcommand \N {\mathcal N}
\renewcommand \O {\mathcal O}
\renewcommand \P {\mathbb P}
\newcommand \W {\mathcal W}
\newcommand \X {\mathcal X}
\newcommand \XX {\mathfrak X}
\newcommand \Z {\mathbb Z}

\newcommand \KN {K\!N}
\newcommand \HH {H\!H}
\DeclareMathOperator \Ext {Ext}
\DeclareMathOperator \cone {cone}
\DeclareMathOperator \RHom {RHom}
\DeclareMathOperator \id {id}
\DeclareMathOperator \Tor {Tor}
\DeclareMathOperator \im {im}
\DeclareMathOperator \supp {supp}

\newcommand \Caldararu {C\u{a}ld\u{a}\-raru}

\newtheorem{thm}{Theorem}[section]
\newtheorem{prop}[thm]{Proposition}
\newtheorem{lem}[thm]{Lemma}
\newtheorem{claim}[thm]{Claim}
\theoremstyle{definition}
\newtheorem{defn}[thm]{Definition}
\newtheorem{rmk}[thm]{Remark}

\numberwithin{equation}{section}
\numberwithin{figure}{section}

\begin{document}
\input intro
\input review
\input sod
\begin{cut} \input windows \end{cut}
\input mukai

\newcommand \httpurl [1] {\href{http://#1}{\nolinkurl{#1}}}
\bibliographystyle{plain}
\bibliography{main}

\end{document}

%% file: intro.tex

\begin{abstract}
Associated to a Mukai flop $X \dashrightarrow X'$ is on the one hand a sequence of equivalences $D^b(X) \to D^b(X')$, due to Kawamata and Namikawa, and on the other hand a sequence of autoequivalences of $D^b(X)$, due to Huybrechts and Thomas.  We work out a complete picture of the relationship between the two.  We do the same for standard flops, relating Bondal and Orlov's derived equivalences to spherical twists, extending a well-known story for the Atiyah flop to higher dimensions.
\end{abstract}

\maketitle

\section*{Introduction}

The Atiyah flop is one of the most-studied objects in derived categories and mirror symmetry.  One favorite fact is the following.  Let $X$ be a complex threefold containing a $(-1,-1)$-curve, that is, a $\P^1$ with normal bundle $\O(-1) \oplus \O(-1)$, and let
\[ \xymatrixresized{
& \tilde X \ar[ld]_q \ar[rd]^p \\
X & & X'
} \]
be its Atiyah flop.  Bondal and Orlov \cite[\S3]{bo} showed that the functor $p_* q^*\colon D^b(X) \to D^b(X')$ is an equivalence.  Symmetrically, $q_* p^*\colon D^b(X') \to D^b(X)$ is an equivalence, but these two natural equivalences are not inverse to one another: their composition $q_* p^* p_* q^*$ is an important autoequivalence of $D^b(X)$, the inverse of the spherical twist around $\O_{\P^1}(-1)$.  One seeks to prove similar ``flop-flop = twist'' results for other classes of flops.

Bondal and Orlov also produced equivalences for flops at $(-2,0)$-curves in threefolds, and Toda \cite{toda} proved a flop-flop = twist result for them by replacing the object $\O_{\P^1}(-1) \in D^b(X)$, which is no longer spherical, with a spherical functor $D^b(\operatorname{Spec}(\C[x]/x^n)) \to D^b(X)$.  Bridgeland \cite{bridgeland_flops} produced an equivalence for general flops of smooth threefolds, encompassing flops at $(-3,1)$-curves and trees of $\P^1$s, and Chen \cite{chen} extended this to certain singular threefolds.  Donovan and Wemyss \cite{dw1,dw2} proved a flop-flop = twist result in this setting which suggests a spherical functor from the derived category of modules over a certain non-commutative algebra.  \begin{cut}We review the definition of spherical and $\P$-functors and the associated autoequivalences (``twists'') in \S\ref{review}.\end{cut}

\subsection*{Standard flops in higher dimensions}
First we investigate the \emph{standard flop}, which generalizes the Atiyah flop to higher dimensions.  Let $X$ be a complex $(2n+1)$-fold containing a $\P^n$ with normal bundle $\O(-1)^{n+1}$, and let $X \xleftarrow{q} \tilde X \xrightarrow{p} X'$ be its flop.  Then $p_* q^*$ and $q_* p^*$ are still equivalences, and $\O_{\P^n}(-1)$ is still a spherical object, but $q_* p^* p_* q^*$ is no longer a spherical twist.  This turns out just to be a 
normalization issue, as follows.  For each $k \in \Z$, define a pair of functors
\[ \xymatrix{
D^b(X) \ar@<1ex>[rrrrrr]^*+{BO_k := p_*\big(\O_{\tilde X}(kE) \otimes q^*(-)\big)} &&&&&&
D^b(X'), \ar@<1ex>[llllll]^*+{BO'_k := q_*\big(\O_{\tilde X}(kE) \otimes p^*(-)\big)}
} \]
where $E \subset \tilde X$ is the exceptional divisor.  Thus in particular $BO_0 = p_* q^*$ and $BO'_0 = q_* p^*$ are the equivalences of Bondal and Orlov, but their proof is easily adapted to show that $BO_k$ and $BO'_k$ are equivalences for all $k$; or see \cite[Prop.~3.1]{kawamata2} for a different proof.  Note that the inverse of $BO_k$ is its left (or right) adjoint, so
\[ BO_k^{-1} = BO'_{n-k}. \]
{ \renewcommand{\thethm}{A}
\begin{thm} \label{BO_single}
Let $X$ be a complex $(2n+1)$-fold containing a $\P^n$ with normal bundle $\O(-1)^{n+1}$, let $X'$ be its flop, and let $BO_k$ and $BO'_k$ be the equivalences defined above.  Then the spherical twist around $\O_{\P^n}(k)$ satisfies
\[ T_{\O_{\P^n}(k)} = BO'_{-k} \circ BO_{n+k+1}, \]
and thus
\[ T_{\O_{\P^n}(k)}^{-1} = BO'_{-k-1} \circ BO_{n+k}. \]
\end{thm} }
\noindent Taking $n=1$ and $k=-1$ in the second statement recovers the ``favorite fact'' that we discussed initially.  For general $n$, we find that
\[ q_* p^* p_* q^* = BO'_0 \circ BO_0 = T_{\O_{\P^n}(-1)}^{-1} \circ T_{\O_{\P^n}(-2)}^{-1} \circ \dotsb \circ T_{\O_{\P^n}(-n)}^{-1}, \]
and indeed that any flop-flop functor $BO'_k \circ BO_l$ can be written as a product of spherical twists or inverse spherical twists around $\O_{\P^n}(m)$ for suitable $m$.

\smallskip
Theorem \ref{BO_single} is a special case of the following family version:
{ \renewcommand{\thethm}{A$'$}
\begin{thm} \label{BO_family}
Suppose we have 
\[ \xymatrixresized{
\P V \ar@{^{(}->}[r]^j \ar[d]_\varpi & X \\
{\phantom,}Z,
} \]
where $Z$ is a smooth complex variety, $V$ is a vector bundle of rank $n+1$ over $Z$, and $j$ is a closed embedding with normal bundle $\N_{\P V/X} = \O_{\P V}(-1) \otimes \varpi^* V'$ for some (possibly different) vector bundle $V'$ of rank $n+1$ over $Z$; thus in particular $\dim X = \dim Z + 2n+1$.  Let $X'$ be the flop of $X$ along $j(\P V)$, and let $BO_k$ and $BO'_k$ be the equivalences defined above.  Then the functor
\[ j_*\big(\O_{\P V}(k) \otimes \varpi^*(-)\big) \colon D^b(Z) \to D^b(X) \]
is spherical, and the associated spherical twist is $BO'_{-k} \circ BO_{n+k+1}$.
\end{thm} }
\begin{cut}
\noindent We give two proofs: one using semi-orthogonal decompositions, close in spirit to Bondal and Orlov's proof that $p_* q^*$ is an equivalence, in \S\ref{sod}, and one using variation of GIT quotients and ``window shifts'' in \S\ref{windows}.
\end{cut}

\subsection*{Mukai flops}
Next we turn our attention to Mukai flops.  Let $X$ be a complex $2n$-fold containing a $\P^n$ whose normal bundle is isomorphic to its cotangent bundle $\Omega^1_{\P^n}$, and let $X \xleftarrow{q} \tilde X \xrightarrow{p} X'$ be its Mukai flop.  If $n>1$ then Kawamata \cite[\S5]{kawamata} and Namikawa \cite{namikawa} showed that $p_* q^*$ is not an equivalence, but it can be modified to give one: the exceptional divisor $E \subset \tilde X$ is naturally identified with the universal hyperplane in $\P^n \times \P^{n*}$, and the correspondence
\[ \hat X := \tilde X \cup_E (\P^n \times \P^{n*}) \]
induces equivalences $D^b(X) \longleftrightarrow D^b(X')$.  Again these are not inverse to one another; but if we replace the bare correspondence $\hat X$ with the line bundle $\L^k$ obtained by gluing together $\O(kE)$ on $\tilde X$ and $\O(-k,-k)$ on $\P^n \times \P^{n*}$, we get equivalences
\[ \xymatrix{
D^b(X) \ar@<1ex>[rrrrrr]^*+{\KN_k := \hat p_*\big(\L^k \otimes \hat q^*(-)\big)} &&&&&&
D^b(X'), \ar@<1ex>[llllll]^*+{\KN'_k := \hat q_*\big(\L^k \otimes \hat p^*(-)\big)}
} \]
where $X \xleftarrow{\hat q} \hat X \xrightarrow{\hat p} X'$ are the obvious maps, and these equivalences satisfy
\[ \KN_k^{-1} = \KN'_{n-k}. \]
{ \renewcommand{\thethm}{B}
\begin{thm} \label{KN_single}
Let $X$ be a smooth complex $2n$-fold containing a $\P^n$ with normal bundle $\Omega^1_{\P^n}$, let $X'$ be its Mukai flop, and let $\KN_k$ and $\KN'_k$ be the equivalences defined above.  Then the $\P$-twist around $\O_{\P^n}(k)$ satisfies
\[ P_{\O_{\P^n}(k)} = \KN'_{-k} \circ \KN_{n+k+1}. \]
\end{thm} }
\noindent Cautis \cite[Prop.~6.8]{cautis_about} proved a special case of this, where $X$ is the total space of $\Omega^1_{\P^n}$ and $k=-n$, as a corollary to an elaborate ``categorical $\mathfrak{sl}_2$-action'' on cotangent bundles of Grassmannians which he and his coauthors had built up over the course of several papers.  Our proof is different: we deduce Theorem \ref{KN_single} from Theorem \ref{BO_single}.
\smallskip

\begin{cut}
The case $n=1$ of Theorem \ref{KN_single} appears to be trivial at first, but in fact it is rather interesting.  In this case the Mukai flop does nothing: $p$ and $q$ are isomorphisms, so $X' = X$.  But the equivalence $\KN_0 = \KN'_0$ is not the identity: it is the inverse of the spherical twist around $\O_{\P^1}(-1)$, and the statement that $\KN'_0 \circ \KN_0$ is the inverse of the $\P$-twist around $\O_{\P^1}(-1)$ amounts to Huybrechts and Thomas's statement \cite[Prop.~2.9]{ht} that the $\P$-twist around a $\P^1$-object is the square of the spherical twist around the same object.  Thus in higher dimensions, while the $\P$-twist around $\O_{\P^n}(-k)$ does not have a literal square root, the Kawamata--Namikawa equivalences to $D^b(X')$ should perhaps be seen as its metaphorical square roots. This is most striking when $n$ is odd, so we can get the same index on $\KN$ and $\KN'$:
\begin{align*}
P_{\O_{\P^n}(-(n+1)/2)} &= \KN'_{(n+1)/2} \circ \KN_{(n+1)/2} \\
P_{\O_{\P^n}(-(n+1)/2)}^{-1} &= \KN'_{(n-1)/2} \circ \KN_{(n-1)/2}.
\end{align*}
One wonders then whether $\P$-twists around other $\P$-objects, e.g.\ the structure sheaf of a hyperk\"ahler variety, can be factored in interesting ways.
\smallskip
\end{cut}

Theorem \ref{KN_single} is a special case of a family version\begin{cut}, which we prove in \S\ref{mukai}\end{cut}:
{ \renewcommand{\thethm}{B$'$}
\begin{thm} \label{KN_family}
Suppose we have 
\[ \xymatrixresized{
\P V \ar@{^{(}->}[r]^j \ar[d]_\varpi & X \\
{\phantom,}Z,
} \]
where $Z$ is a smooth complex projective variety with $\HH^\text{odd}(Z) = 0$,\begin{cut}\footnote{See Remark \ref{rmk_on_HH_odd} for a discussion of this assumption.}\end{cut} $V$ is a vector bundle of rank $n+1$ over $Z$, and $j$ is a closed embedding with normal bundle $\N_{\P V/X} = \Omega^1_{\P V/Z}$; thus in particular $\dim X = \dim Z + 2n$.  Let $X'$ be the Mukai flop of $X$ along $j(\P V)$, and let $\KN_k$ and $\KN'_k$ be the equivalences defined above.  Then the functor
\[ j_*\big(\O_{\P V}(k) \otimes \varpi^*(-)\big) \colon D^b(Z) \to D^b(X) \]
is a $\P^n$-functor, and the associated $\P$-twist is $\KN'_{-k} \circ \KN_{n+k+1}$.
\end{thm} }

\begin{cut}
Theorems \ref{BO_family} and \ref{KN_family} can be generalized further to allow a $\P^n$-bundle over $Z$ that is not the projectivization of a vector bundle, but we omit this generalization for the sake of clarity and for want of a compelling example.
\end{cut}

In \cite{torsion_sheaves} we apply Theorem \ref{KN_family} to the following example: $Z$ is a general K3 surface of degree $2n$, $\P V$ is the total space of the universal hyperplane section of $Z$, which is both a $\P^n$-bundle over $Z$ and a family of genus-$(n+1)$ curves over $\P^{n+1}$, and $j$ is the Abel--Jacobi embedding into the relative Jacobian of the latter family.

\begin{cut}
It would be very interesting to find something like a spherical or $\P$-functor associated to Markman's stratified Mukai flops \cite{markman_brill_noether}, and to relate the corresponding twist to the equivalences of Cautis, Kamnitzer, and Licata \cite{ckl, cautis_stratified}.
\end{cut}

\subsection*{Conventions} We work throughout with smooth quasi-projective varieties over $\C$; these hypotheses can certainly be relaxed, but we leave that task to the interested reader.  All pushforwards, tensor products, etc.\ are implicitly derived.  We freely identify Fourier--Mukai functors $D^b(X) \to D^b(Y)$ with their kernels in $D^b(X \times Y)$, and units and counits of adjunctions with certain natural maps between kernels, so we are implicitly working in some version of the 2-category $\mathcal Var$ of \Caldararu\ and Willerton \cite{andrei_simon}.

\begin{cut}
If $V$ is a vector bundle over $Z$ then $\varpi\colon \P V \to Z$ is the projective bundle of 1-dimensional \emph{subspaces} of the fibers of $V$, and $\O_{\P V}(-1)$ is the tautological sub-bundle of $\varpi^* V$.  In particular, if $L$ is a line bundle then $\P(V \otimes L) = \P V$ but $\O_{\P(V \otimes L)}(-1) = \O_{\P V}(-1) \otimes \varpi^* L$.
\end{cut}

\subsection*{Acknowledgements}  We thank Daniel Huybrechts, Richard Thomas, and Arend Bayer for helpful discussions, and Sabin Cautis for expert advice on \cite{cautis_about}.  Much of this work was done while visiting the Hausdorff Research Institute for Mathematics in Bonn, Germany, and we thank them for their hospitality.  N.A.\ was partially supported by NSF grant no.\ DMS--0905923.  W.D.\ was supported by World Premier International Research Center Initiative (WPI Initiative), MEXT, Japan; and by EPSRC grant~EP/G007632/1.  C.M.\ was supported by the EPSRC Doctoral Prize Research Fellowship Grant EP/K503034/1.

%% file: review.tex

\section{Review of spherical and \texorpdfstring{$\P$}{P}-twists} 
\label{review}

\subsection*{Spherical twists}
An object $\E \in D^b(X)$ is called \emph{spherical} if we have
\[ \Ext^*_X(\E,\E) = H^*(S^n, \C) = \begin{cases}
\C & i = 0 \text{ or } n \\
0 & \text{otherwise,}
\end{cases} \]
where $n = \dim X$, and $\E \otimes \omega_X \cong \E$.  \begin{cut}The main example for us is the sheaf $\O_{\P^n}(k)$ appearing in Theorem \ref{BO_single}, but line bundles and other rigid stable vector bundles on K3 surfaces provide another important class of examples.  \end{cut} Seidel and Thomas \cite{st} showed that if $\E$ is spherical then the \emph{spherical twist} $T_\E \colon D^b(X) \to D^b(X)$ given by
\[ T_\E(\F) = \cone(\E \otimes \RHom(\E,\F) \xrightarrow{\text{ eval }} \F) \]
is an equivalence.
\begin{ins}
A relative version of this notion is as follows:
\end{ins}
\begin{cut}
Spherical twists arise naturally in mathematical string theory, as monodromy operators associated to loops in the complexified K\"ahler moduli space around limit points where a subvariety $Y \subset X$ is contracted to a point; the pushforward of a line bundle on $Y$ is often a spherical object in $D^b(X)$.  When $Y$ is instead contracted to a positive-dimensional variety $Z$, we need a more general construction, Horja's EZ-spherical twist \cite{horja}.  This is most conveniently described in the even more general language of spherical functors, which we now quickly review; for a more detailed discussion we suggest \cite[\S1]{nick}.
\end{cut}

\begin{defn}
Let $Z$ and $X$ be smooth, quasi-projective varieties,\footnote{This is not the most general possible setting, but it is sufficient for our purposes here.} let $F\colon D^b(Z) \to D^b(X)$ be a Fourier--Mukai functor induced by a kernel whose support is proper over $Z$ and $X$, and let $L, R\colon D^b(X) \to D^b(Z)$ be the left and right adjoints of $F$.  Define the \emph{cotwist} $C\colon D^b(Z) \to D^b(Z)$ and the \emph{twist} $T\colon D^b(X) \to D^b(X)$ as the cones on the unit and counit of the adjunction:
\begin{align*}
C &:= \cone(\id_Z \xrightarrow\eta RF) &
T &:= \cone(FR \xrightarrow\epsilon \id_X).
\end{align*}
Then $F$ is called \emph{spherical} if $C$ is an equivalence and $R \cong CL$.
\end{defn}

A spherical object is then the same as a spherical functor from $D^b$(point).

\begin{thm}[{Rouquier \cite{rouquier}, Anno \cite{anno}, Anno and Logvinenko \cite{al}; see also \cite[Thm.~1]{nick}}] 
If $F$ is spherical then $T$ is an equivalence.
\end{thm}

\begin{thm}[Horja \cite{horja}] \label{horja_thm}
\begin{cut}Suppose we have \end{cut}\begin{ins}Take smooth quasi-projective varieties\end{ins}
\[ \xymatrixresized{
Y \ar@{^{(}->}[r]^j \ar[d]_\varpi & X \\
{\phantom,}Z,
} \]
where \begin{cut}$X$, $Y$, and $Z$ are smooth quasi-projective varieties, \end{cut}$j$ is a closed embedding of codimension $d>0$\begin{cut} with normal bundle $\N_{Y/X}$\end{cut}, and $\varpi$ is proper.  Suppose \begin{cut}that an object \end{cut}$\E \in D^b(Y)$ satisfies
\[ \varpi_*(\E^* \otimes \E \otimes \Lambda^i \N_{Y/X}) = \begin{cases}
\O_Z & i = 0 \\
0 & 0 < i < d,
\end{cases} \]
and $\E \otimes j^* \omega_X \cong \E \otimes \varpi^* \L$ for some line bundle $\L$ on $Z$.
Then the functor
\[ F := j_*\big(\E \otimes \varpi^*(-)\big)\colon D^b(Z) \to D^b(X) \]
is spherical\begin{cut} with cotwist $C = - \otimes \L^* \otimes \omega_Z[\dim Z - \dim X]$\end{cut}.
\end{thm}
\noindent The hypotheses imply that
\[ \varpi_*(\E^* \otimes \E \otimes \Lambda^d \N_{Y/X}) = \L^* \otimes \omega_Z [\dim Z - \dim Y], \]
using Grothendieck duality.  Note that the twist $T\colon D^b(X) \to D^b(X)$ acts as the identity on objects supported on $X \setminus j(Y)$.  Horja also considered the case $d = 0$, but this is not relevant for us.

Horja's Theorem \ref{horja_thm} in fact predates the language of spherical functors.  He calls an object \emph{EZ-spherical}\footnote{To explain the name, we mention that Horja used $E$ for the space we are calling $Y$.  We want to reserve $E$ for an exceptional divisor later.} if it satisfies the hypotheses of the theorem,  Huybrechts \cite[Def.~8.43]{huybrechts_fm} calls an object EZ-spherical if the functor $F$ is spherical, although he does not use this language either.  \begin{cut}He claims [\emph{ibid}., Rmk.~8.50] that this definition is slightly weaker than Horja's, but of course it is much harder to check in practice.\end{cut}

\begin{cut}
In the present paper we only use spherical functors of Horja's form, but as we mentioned in the introduction, other classes of flops lead one to consider spherical functors where $Z$ is replaced with a fat point, or with a non-commutative algebra; and examples of a very different flavor arise in studying hyperk\"ahler 4-folds \cite{nick,ciaran_kummer} or canonical covers of Hilbert schemes of Enriques surfaces \cite{krug_sosna}.
\end{cut}
\smallskip

We give a short proof of Theorem \ref{horja_thm} when $\E$ is a line bundle, which will serve as a model of our proof of Proposition \ref{P-prop}.  The main ideas are the same as in Horja's original proof, but our assumption that $\E$ is a line bundle yields many simplifications, as does the spherical functor language.  We begin by stating a useful lemma, which generalizes \cite[Cor.~11.4(ii) and Prop.~11.8]{huybrechts_fm}:

\begin{lem} \label{j*j*}
Let $j \colon Y \hookrightarrow X$ be a closed embedding of smooth quasi-projective varieties, with normal bundle $\N$.  Then $j^* j_*$ (resp.\ $j^! j_*$) is given by a kernel with cohomology sheaves $\H^i = \Delta_* \Lambda^{-i} \N^*$ (resp.\ $\H^i = \Delta_* \Lambda^i \N$).
\begin{ins} \qed \end{ins}
\end{lem}

\begin{cut}
As a consequence, for a sheaf $\F$ on $Y$ we get $\H^i(j^* j_* \F) = \F \otimes \Lambda^{-i} \N^*$ and $\H^i(j^! j_* \F) = \F \otimes \Lambda^i \N$, and for a general object $\F \in D^b(Y)$ we get a spectral sequence.  Arinkin and \Caldararu\ \cite{ac} have shown that if the normal bundle sequence $0 \to TY \to TX|_Y \to \N \to 0$ splits, then the kernels inducing $j^*j_*$ and $j^! j_*$ are formal, i.e.\ they split as the sum of their cohomology sheaves.

\begin{proof}[Proof of Lemma \ref{j*j*}]
Let $\Gamma \subset Y \times X$ be the graph of $j$, and let $\bar\Gamma \subset X \times Y$ be its transpose.  Then the kernel inducing $j^* j_*$ is obtained by taking
\[ \O_{\Gamma \times Y} \otimes \O_{Y \times \bar\Gamma} \in D^b(Y \times X \times Y) \]
and pushing down to $Y \times Y$.  The intersection
\[ (\Gamma \times Y) \cap (Y \times \bar\Gamma) \]
is the image of
\[ (1 \times j \times 1)\colon Y \to Y \times X \times Y. \]
It is not a transverse intersection, but it is smooth, and we find that the excess normal bundle is exactly $\N$.  Thus we have
\begin{equation} \label{tors}
\Tor_i^{Y \times X \times Y}(\O_{\Gamma \times Y}, \O_{Y \times \bar\Gamma}) = (1 \times j \times 1)_* \Lambda^i \N^*;
\end{equation}
see for example \cite[VII, Prop.~2.5]{SGA6} or \cite[Prop.~A.5]{cks}.  Pushing the sheaves \eqref{tors} down to $Y \times Y$ we get $\Delta_* \Lambda^i \N$, so the claim about $j^* j_*$ follows using the Grothendieck spectral sequence in the form \cite[Eq.~2.6]{huybrechts_fm}.

The claim about $j^! j_*$ follows from the fact that
\[ j^! j_* = j^* j_* \otimes \det \N[\dim Y - \dim X]. \qedhere \]
\end{proof}
\end{cut}

\begin{proof}[Proof of Theorem \ref{horja_thm} when $\E$ is a line bundle] \label{horja_pf_with_lb}
We have
\begin{align*}
F &= j_*\big(\E \otimes \varpi^*(-)\big) &
\text{and} & &
R &= \varpi_*\big(\E^* \otimes j^!(-)\big).
\end{align*}
We will show that the kernel inducing $R F$ has cohomology sheaves
\begin{equation} \label{RF_claim}
\H^i = \begin{cases}
\Delta_* \O_Z & i = 0 \\
\Delta_*(\L^* \otimes \omega_Z) & i = \dim X - \dim Z \\
0 & \text{otherwise.}
\end{cases} \end{equation}
It follows that the cotwist $C = - \otimes \L^* \otimes \omega_Z[\dim Z - \dim X]$, which is an equivalence, and the verification that $R \cong C L$ is straightforward.

By Lemma \ref{j*j*}, the kernel inducing $j^! j_*$ has cohomology sheaves
\begin{equation} \label{Hjj}
\H^i = \Delta_* (\Lambda^i \N) \in \operatorname{Coh}(Y \times Y).
\end{equation}
To get the kernel inducing $\E^* \otimes j^! j_*(\E \otimes -)$ we tensor with $\E \boxtimes \E^*$, which does not change the sheaves \eqref{Hjj} because $\E$ is a line bundle.  To get the kernel inducing $R F$ we apply $(\varpi \times \varpi)_*$, which sends the sheaves \eqref{Hjj} to $\Delta_* \O_Z$ for $i=0$ and to zero for $0 < i < d$ by hypothesis, and to $\Delta_*(\L^* \otimes \omega_Z)[\dim Z - \dim Y]$ for $i=d$ as we remarked earlier.  Thus we get \eqref{RF_claim} using the Grothendieck spectral sequence\begin{cut}, again\end{cut} in the form \cite[Eq.~2.6]{huybrechts_fm}.
\end{proof}
\bigskip

\subsection*{\texorpdfstring{$\P$}{P}-twists}
An object $\E \in D^b(X)$ is called a \emph{$\P$-object} if we have
\[ \Ext^*_X(\E,\E) = H^*(\P^n, \C) = \C[h]/h^{n+1}, \quad \deg h = 2, \]
as rings, where $\dim X = 2n$, and $\E \otimes \omega_X \cong \E$.  \begin{cut}The main example for us is the sheaf $\O_{\P^n}(k)$ appearing in Theorem \ref{KN_single}; the other main example is a line bundle on a hyperk\"ahler $2n$-fold.  \end{cut}Huybrechts and Thomas \cite{ht} showed that if $\E$ is a $\P$-object then the \emph{$\P$-twist} $P_\E\colon D^b(X) \to D^b(X)$ given by a certain double cone
\begin{equation} \label{ht_double_cone}
P_\E(\F) = \cone\!\big(\cone\!\big(\E \otimes \RHom(\E,\F)[-2] \to \E \otimes \RHom(\E,\F)\big) \xrightarrow{\text{eval}} \F\;\big)
\end{equation}
is an equivalence.

Later we will need the following instance of \cite[Prop.~1.4]{ht}, which relates $\P$-twists to spherical twists:
\begin{prop}[Huybrechts, Thomas] \label{intertwinement_prop}
Let $X$ be the total space of $\Omega^1_{\P^n}$, let $\X$ be the total space of $\O_{\P^n}(-1)^{n+1}$, and let $\iota\colon X \to \X$ be the embedding given by the Euler sequence.  Then the object $\O_{\P^n}(k) \in D^b(X)$, supported on the zero section, is a $\P$-object; its pushforward $\iota_* \O_{\P^n}(k) \in D^b(\X)$ is a spherical object; and the spherical and $\P$-twists satisfy
\[ T_{\iota_* \O_{\P^n}(k)} \circ \iota_* \cong \iota_* \circ P_{\O_{\P^n}(k)}. \]
\end{prop}

\begin{cut}
\newpage
\end{cut}
In \cite{nick}, the first author proposed a definition of \emph{$\P$-functors}:
\begin{defn} \label{P-def}
Let $Z$ and $X$ be smooth, quasi-projective varieties,\footnote{Again this is not the most general possible setting.} let $F\colon D^b(Z) \to D^b(X)$ be a Fourier--Mukai functor induced by a kernel whose support is proper over $Z$ and $X$, and let $L, R\colon D^b(X) \to D^b(Z)$ be the left and right adjoints of $F$.  Then $F$ is called a \emph{$\P$-functor} if there is an autoequivalence $H\colon D^b(Z) \to D^b(Z)$ such that
\[ RF \cong \id_Z \oplus H \oplus H^2 \oplus \dotsb \oplus H^n, \]
and \begin{cut}the monad structure $RFRF \xrightarrow{R\epsilon F} RF$ looks like the multiplication in $H^*(\P^n,\C)$, in the following sense: \end{cut}the composition
\begin{cut}\[ HRF \hookrightarrow RFRF \xrightarrow{R\epsilon F} RF, \]\end{cut}
\begin{ins}$ HRF \hookrightarrow RFRF \xrightarrow{R\epsilon F} RF, $\end{ins}
when written in components
\[ H \oplus H^2 \oplus \dotsb \oplus H^n \oplus H^{n+1} \to \id_Z \oplus H \oplus H^2 \oplus \dotsb \oplus H^n, \]
is of the form
\[ \begin{pmatrix}
* & * & \cdots & * & * \\
1 & * & \cdots & * & * \\
0 & 1 & \cdots & * & * \\
\vdots & \vdots & \ddots & \vdots & \vdots \\
0 & 0 & \cdots & 1 & *
\end{pmatrix}; \]
and moreover $R \cong H^n L$.
\end{defn}
He then constructed an equivalence $P\colon D^b(X) \to D^b(X)$, as a certain double cone
\begin{equation} \label{my_double_cone}
P = \cone\!\big(\cone(FHR \to FR) \xrightarrow\epsilon \id_X\big).
\end{equation}
This recovers Huybrechts and Thomas's equivalence if we take $Z$ to be a point and $H = [-2]$.  \begin{cut}For Huybrechts and Thomas there is a unique way to take the double cone \eqref{ht_double_cone}, but in general there is a somewhat subtle choice to be made.

The example of Theorem \ref{KN_family} and Krug's examples \cite{krug_diagonal, krug_nakajima} are ``Horja-esque'' in that the resulting $\P$-twists act as the identity on a Zariski open subset; the examples of Addington \cite{nick} and Meachan \cite{ciaran_kummer}, in contrast, do not act as the identity on any skyscraper sheaf.

\begin{rmk} \label{P-rmk}
Definition \ref{P-def} is probably not quite the right one.  Rather than requiring $RF$ to split, it would be better to require a filtration of $RF$ with quotients $\id_Z, H, H^2, \dotsc, H^n$; but then it is hard to make precise the hypothesis that ``the monad structure $RFRF \xrightarrow{R\epsilon F} RF$ looks like $H^*(\P^n,\C)$.''  Cautis's definition of a $\P$-functor \cite[\S6.2]{cautis_about} solves this problem at the cost of allowing only $H = [-2]$.  Maybe this is not such a great cost: in all known examples we have $H = - \otimes \L[-k]$ for some line bundle $\L$ on $Z$ and some $k > 0$, and Cautis's approach extends easily to this case.  The issue of how to take the double cone \eqref{my_double_cone} remains subtle, however, and Cautis has to make the somewhat artificial assumption that $\HH^1(Z) = 0$; see \cite[Rmk.~6.7]{cautis_about}.

On the other hand, E.~Segal \cite[\S4]{ed_spherical} has observed that a $\P$-object is the same as a spherical functor from the category of dg-modules over the dg-algebra $\C[h]$, where $h$ has homological degree 2.\footnote{By Koszul duality, this is equivalent to the category of dg-modules over $\C[\epsilon]/\epsilon^2$, where $\epsilon$ has homological degree $-1$; thus Segal's observation agrees with \cite[\S3.2, Example~6]{nick}.}  Ultimately, the way forward with $\P$-functors will probably be to work out a family version of Segal's observation that covers the known examples, and then retire them in favor of spherical functors from dg bases; but we are not ready to take this step yet.
\end{rmk}
\end{cut}

%% file: sod.tex

\begin{ins}
\section{Standard flops} \label{sod}
\end{ins}
\begin{cut}
\section{Standard flops via semi-orthogonal decompositions} \label{sod}
\end{cut}

Assume the set-up of Theorem \ref{BO_family}: we have
\[ \xymatrixresized{
\P V \ar@{^{(}->}[r]^j \ar[d]_\varpi & X \\
{\phantom,}Z,
} \]
and $\N_{\P V/X} = \O_{\P V}(-1) \otimes \varpi^* V'$.  Thus we see that for any $k \in \Z$, the line bundle $\O_{\P V}(k)$ satisfies the hypotheses of Theorem \ref{horja_thm} with $\L = \det V^* \otimes \det V'^* \otimes \omega_Z$, so the functor
\[ F_k := j_*\big(\O_{\P V}(k) \otimes \varpi^*(-)\big)\colon D^b(Z) \to D^b(X) \]
is spherical.  In this section we will show that the associated spherical twist $T_k$ satisfies
\[ T_k \circ BO'_{-k-1} \cong BO'_{-k}. \]
Since $(BO'_{-k-1})^{-1} = BO_{n+k+1}$ this implies Theorem \ref{BO_family}.
\smallskip

Following Huybrechts \cite[\S11.3]{huybrechts_fm}, we adopt the notation
\[ \xymatrixresized{
& & \tilde X \ar[lldd]_q \ar[rrdd]^p \\
& & E \ar[ld]_{\pi} \ar@{_(->}[u]^i \ar[rd]^{\pi'} \\
X & \P V \ar@{_(->}[l]^j \ar[rd]_\varpi & & \P V' \ar@{^(->}[r]_{j'} \ar[ld]^{\varpi'} & X' \\
& & {\phantom.}Z.
} \]
The exceptional divisor $E \subset \tilde X$ is identified with $\P V \times_Z \P V'$, and under this identification we have $\O_E(E) = \O_E(-1,-1)$.  For brevity we define
\[ M_m := \O_{\tilde X}(m E) \otimes - \colon D^b(\tilde X) \to D^b(\tilde X), \]
so for example $BO'_{-k} = q_* M_{-k} p^*$.

The basic observation we wish to exploit is the following.  Suppose for a moment that $Z$ is a point.   Then the spherical object $j_* \O_{\P V}(k) \in D^b(X)$ is $q_*$ of the exceptional object $i_* \O_E(k,0) \in D^b(\tilde X)$, and the formula for the spherical twist around the former is very similar to the formula for the mutation past the latter.  This motivation gets a bit lost in the proof below, but it might be glimpsed in the proof of Claim \ref{second_claim}.

When $Z$ is general we replace the exceptional objects $i_* \O_E(a,b)$, where $a,b\in\Z$, with the images of the fully faithful functors
\[ I_{a,b} := i_* \big(\O_E(a,b) \otimes \pi'^* \varpi'^*(-)\big)\colon D^b(Z) \to D^b(\tilde X). \]
We make a few preparatory observations about $I_{a,b}$.  First,
\[ \setlength \arraycolsep {1.6pt}
\begin{cut}\renewcommand \arraystretch {1.2}\end{cut}
\begin{array}{rlp{2cm}lrl}
p_* I_{a,b} &= 0 & & \text{for} & a &= -n, \dotsc,-1 \\
q_* I_{a,b} &= 0 & & \text{for} & b &= -n, \dotsc,-1,
\end{array} \]
and with a little more work,
\[ \setlength \arraycolsep {1.6pt}
\begin{cut}\renewcommand \arraystretch {1.2}\end{cut}
\begin{array}{rlp{2cm}lrl}
I_{a,b}^R p^* &= 0 & & \text{for} & a &= 0, \dotsc,n-1 \\
I_{a,b}^R q^* &= 0 & & \text{for} & b &= 0, \dotsc,n-1,
\end{array} \]
where $I_{a,b}^R$ is the right adjoint of $I_{a,b}$.  Next, let $L_{a,b}\colon D^b(\tilde X) \to D^b(\tilde X)$ be the left mutation past $\im(I_{a,b}) \subset D^b(\tilde X)$:
\[ L_{a,b} := \cone(I_{a,b} I_{a,b}^R \xrightarrow{\epsilon} 1). \]
Then we have
\[ \setlength \arraycolsep {1.6pt}
\begin{cut}\renewcommand \arraystretch {1.2}\end{cut}
\begin{array}{rlp{2cm}lrl}
p_* L_{a,b} &= p_* & & \text{for} & a &= -n, \dotsc,-1 \\
q_* L_{a,b} &= q_* & & \text{for} & b &= -n, \dotsc,-1, \\
L_{a,b} p^* &= p^* & & \text{for} & a &= 0, \dotsc,n-1 \\
L_{a,b} q^* &= q^* & & \text{for} & b &= 0, \dotsc,n-1.
\end{array} \]
Finally we observe that $M_m I_{a,b} = I_{a-m,b-m}$, and thus
\[ M_m L_{a,b} = L_{a-m,b-m} M_m. \]

Now take the exact sequence
\[ 0 \to \O_{\tilde X}(-E) \to \O_{\tilde X} \to i_* \O_E \to 0, \]
tensor with $p^*(-)$, and manipulate the third term to get an exact triangle of functors
\[ M_{-1} p^* \to p^* \to i_* \pi'^* j'^*. \]
We claim that applying $q_* M_{-k} L_{0,-k}$ annihilates the third term, so we are left with an isomorphism
\begin{equation} \label{iso_to_simplify}
q_* M_{-k} L_{0,-k} M_{-1} p^* \xrightarrow\cong q_* M_{-k} L_{0,-k} p^*.
\end{equation}
Accepting the claim for a moment, we can simplify \eqref{iso_to_simplify} using the above preparations to get
\[ q_* L_{k,0} M_{-k-1} p^* \xrightarrow\cong q_* M_{-k} p^*. \]
The right-hand side is $BO'_{-k}$, and we claim the left-hand side is $T_k BO'_{-k-1}$.  Once we prove these two claims we will have proved the theorem.
\smallskip 

\begin{claim}
$q_* M_{-k} L_{0,-k} \circ i_* \pi'^* = 0$.
\end{claim}
\begin{proof}
From the Beilinson semi-orthogonal decomposition of $D^b(\P V')$ and the fact that $i_* \pi'^*$ is fully faithful we get a semi-orthogonal decomposition
\[ \im(i_* \pi'^*) = \bigl \langle \ \im(I_{0,-k-n}), \ \dotsc, \ \im(I_{0,-k-1}), \ \im (I_{0,-k}) \ \bigr\rangle. \]
Now $L_{0,-k}$ annihilates the last factor and acts as the identity on the others, and $q_* M_{-k}$ annihilates the other factors.
\end{proof}

\begin{claim} \label{second_claim}
There is a natural isomorphism of functors
\begin{equation} \label{second_claim_eq}
q_* L_{k,0} M_{-k-1} p^* \cong T_k BO'_{-k-1}.
\end{equation}
\end{claim}
\begin{proof}
Observe that
\[ F_k = q_* I_{k,0}. \]
Thus the claim is that the left-hand side of \eqref{second_claim_eq}, which is
\begin{equation} \label{first_cone}
\cone( q_* I_{k,0} I_{k,0}^R M_{-k-1} p^* \xrightarrow{q \epsilon M p} q_* M_{-k-1} p^* ),
\end{equation}
is isomorphic to the right-hand side, which is
\begin{equation} \label{second_cone}
\cone(q_* I_{k,0} I_{k,0}^R q^! q_* M_{-k-1} p^* \xrightarrow{\epsilon q M p} q_* M_{-k-1} p^*).
\end{equation}
It is enough to show that the map
\begin{equation} \label{crucial_iso}
I_{k,0}^R M_{-k-1} p^* \xrightarrow{I \eta M p} I_{k,0}^R q^! q_* M_{-k-1} p^*
\end{equation}
is an isomorphism: apply the octahedral axiom to the triangle
\[ \xymatrix{
q_* I_{k,0} I_{k,0}^R M_{-k-1} p^*
  \ar[d]_{q I I \eta M p}
  \ar[rrd]^{q \epsilon M p} \\
q_* I_{k,0} I_{k,0}^R q^! q_* M_{-k-1} p^*
  \ar[rr]_{\epsilon q M p}
& & q_* M_{-k-1} p^*, } \]
where the diagonal arrow is the one in \eqref{first_cone}, the horizontal arrow is the one in \eqref{second_cone}, and the vertical arrow is $q_* I_{k,0}$ applied to our alleged isomorphism \eqref{crucial_iso}; the triangle commutes because of
\[ \xymatrix{
q_* I_{k,0} I_{k,0}^R
  \ar[d]_{q I I \eta}
  \ar[rr]^{q \epsilon}
& & q_* 
  \ar[d]_{q \eta}
  \ar@{=}[rrd] \\
q_* I_{k,0} I_{k,0}^R q^! q_*
  \ar[rr]_{q \epsilon q q}
& & q_* q^! q_*
  \ar[rr]_{\epsilon q}
& & q_*. } \]

So the claim is that \eqref{crucial_iso} is an isomorphism, or equivalently (taking left adjoints) that the map
\[ p_! M_{k+1} q^* q_* I_{k,0} \xrightarrow{p M \epsilon I} p_! M_{k+1} I_{k,0} \]
is an isomorphism, or equivalently, that
\begin{equation} \label{thing_to_kill}
\cone(q^* q_* I_{k,0} \xrightarrow{\epsilon I} I_{k,0})
\end{equation}
becomes zero when we compose with $p_! M_{k+1}$ on the left.

For clarity we now take $n=3$, but the generalization to arbitrary $n$ is straightforward.  Take the standard semi-orthogonal decomposition for a blow-up
\begin{align*}
D^b(\tilde X) = \langle &\im(M_{-3} i_* \pi^*), \\
                        &\im(M_{-2} i_* \pi^*), \\
                        &\im(M_{-1} i_* \pi^*), \\
                        &\im(q^*) \rangle,
\end{align*}
and expand it, using Beilinson's semi-orthogonal decomposition of $D^b(\P V)$, as follows:
\[ \setlength \arraycolsep {1pt}
\begin{array}{lllllll}
\langle & \im(I_{k-2,-3}), & \im(I_{k-1,-3}), & \im(I_{k,-3}), & \im(I_{k+1,-3}) \\
& & \im(I_{k-1,-2}), & \im(I_{k,-2}), & \im(I_{k+1,-2}), & \im(I_{k+2,-2}) \\
& & & \im(I_{k,-1}), & \im(I_{k+1,-1}), & \im(I_{k+2,-1}), & \im(I_{k+3,-1}), \\
& \im(q^*) \ \rangle.
\end{array} \]
We can reorder this as
\[ \begin{blockarray}{lllll}
\begin{block}{llll\}l}
\langle & \im(I_{k-2,-3}), & \im(I_{k-1,-3}), & \im(I_{k,-3}), \\
& & \im(I_{k-1,-2}), & \im(I_{k,-2}), & \mathcal A \\
& & & \im(I_{k,-1}), \\
\end{block}
\begin{block}{llll\}l}
& \im(I_{k+1,-3}) \\
& \im(I_{k+1,-2}), & \im(I_{k+2,-2}) & & \mathcal B \\
& \im(I_{k+1,-1}), & \im(I_{k+2,-1}), & \im(I_{k+3,-1}), \\
\end{block}
& \im(q^*) \ \ \rangle
\end{blockarray} \]
thanks to the appropriate $\im(I_{a,b})$s being both left and right orthogonal to one another.  We find that $\im(I_{k,0}) \subset {}^\perp \mathcal A$, so the image of \eqref{thing_to_kill} is contained in ${}^\perp \mathcal A$.  Moreover the image of \eqref{thing_to_kill} is the right mutation of $\im(I_{k,0})$ past $\im(q^*)$, so it is contained in $\im(q^*)^\perp$.  Thus it is contained in $\mathcal B$; but $p_! M_{k+1} = p_* M_{n+k+1}$ annihilates $\mathcal B$, which completes the proof.
\end{proof}

\begin{cut}
The end of the proof suggests that one should be able to fit this into the framework of Halpern-Leistner and Shipman's \cite[Thm.~3.11]{dhl_shipman}, but we could not make this work for $n > 1$.
\end{cut}

%% file: windows.tex

\section{Standard flops via window shifts} \label{windows}

In this section we give an alternative proof of Theorem \ref{BO_family} using variation-of-GIT methods developed by E.~Segal and 
the second author \cite{ed_windows,ds12,ds13}, and in much greater generality by Halpern-Leistner and Shipman \cite{dhl,dhl_shipman} and Ballard, Favero, and Katzarkov \cite{bfk}, inspired by a physics analysis due to Herbst, Hori, and Page \cite{hhp}.

We work in the ``local model'' where $X$ is the total space of $\O_{\P V}(-1) \otimes \varpi^* V'$, so $X'$ is the total space of $\O_{\P V'}(-1) \otimes \varpi'^* V$ and the flop $X \dashrightarrow X'$ can be realized as a variation of GIT quotient, as follows.  Let $\C^*$ act on the total space of the vector bundle $V \oplus V'$ over $Z$, with weight 1 on the fibers of $V$ and weight $-1$ on the fibers of $V'$:
\[ \lambda \cdot (v,v') = (\lambda v, \lambda^{-1} v'). \]
Then $X$ and $X'$ are the two GIT quotients of $V \oplus V'$, obtained by removing the unstable loci $0 \oplus V'$ and $V \oplus 0$.  The Artin stack
\[ \XX := \left[\ V \oplus V'\ /\ \C^*\ \right] \]
contains $X$ and $X'$ as open substacks.  We will describe for each $k \in \Z$ a ``window'' subcategory $\W_k \subset D^b(\XX)$ such that the restriction functors $\W_k \to D^b(X)$ and $\W_k \to D^b(X')$ are equivalences.  Next we will show that the resulting equivalence $\psi_k\colon D^b(X) \to D^b(X')$ coincides with the Bondal--Orlov equivalence $BO_{n+k}$.  Finally we will show that the ``window shift'' autoequivalence $\psi_k^{-1} \circ \psi_{k+1}$ of $D^b(X)$ coincides with the spherical twist of Theorem \ref{BO_family}, which completes the proof.  The result for general $X$ can be deduced from the local model by deformation to the normal bundle, as we do for Mukai flops in Proposition \ref{defo_to_normal_cone} in the next section.

\subsection*{More details of the local model}  If $X$ is the total space of $\O_{\P V}(-1) \otimes \varpi^* V'$ then $X'$ is the total space of $\O_{\P V'}(-1) \otimes \varpi'^* V$ and $\tilde X$ is the total space of $\O_E(-1,-1) = \O_{\P V \times_Z \P V'}(-1,-1)$, and we can add projections $r$, $r'$, $\rho$ to the big diagram from earlier:
\[ \xymatrix{
& & \tilde X \ar[lldd]_q \ar[rrdd]^p \ar@<-.5ex>[d]_\rho \\
& & E \ar[ld]_{\pi} \ar@<-.5ex>@{_(->}[u]_i \ar[rd]^{\pi'} \\
X \ar@<-.5ex>[r]_r & \P V \ar@<-.5ex>@{_(->}[l]_j \ar[rd]_\varpi & & \P V' \ar@<.5ex>@{^(->}[r]^{j'} \ar[ld]^{\varpi'} & X' \ar@<.5ex>[l]^{r'} \\
& & {\phantom.}Z.
} \]
Thus the line bundles $\O_{\P V}(1)$, $\O_{\P V'}(1)$, and $\O_E(1,1)$ have natural extensions to $X$, $X'$, and $\tilde X$:
\begin{align*}
\O_X(1) &:= r^* \O_{\P V}(1) &
\O_{X'}(1) &:= r'^* \O_{\P V'}(1) &
\O_{\tilde X}(1,1) &:= \rho^* \O_E(1,1)
\end{align*}
Moreover an integer $l$ determines a character of $\C^*$ and thus a line bundle on $\XX$, and we have $\O_\XX(l)|_X = \O_X(l)$, but $\O_\XX(l)|_{X'} = \O_{X'}(-l)$.

Another important feature of the local model is that the equivalences $BO_k\colon D^b(X) \to D^b(X')$ are $D^b(Z)$-linear: that is, for $\F \in D^b(X)$ and $\G \in D^b(Z)$ we have a functorial isomorphism
\[ BO_k(\F \otimes r^* \varpi^* \G) \cong BO_k(\F) \otimes r'^* \varpi'^* \G. \]

\subsection*{Windows and window equivalences}
For $k \in \Z$, we define the \emph{window} subcategory
\[ \W_k \subset D^b(\XX) \]
to be the full subcategory of $D^b(\XX)$ split-generated over $D^b(Z)$ by
\begin{equation} \label{window_generators}
\O_\XX(k), \O_\XX(k+1), \dotsc, \O_\XX(k+n),
\end{equation}
that is, split-generated by objects of the form $\varrho^* \G \otimes \O_\XX(l)$, where $\varrho\colon \XX \to Z$ is the projection, $\G \in D^b(Z)$, and $k \le l \le k+n$.

If $\iota\colon X \to \XX$ and $\iota'\colon X' \to \XX$ denote the two inclusions then $\iota^*\colon \W_k \to D^b(X)$ and $\iota'^*\colon \W_k \to D^b(X')$ are equivalences.  This follows rapidly from the fact that the direct sum of the generators \eqref{window_generators} of $\W_k$ restricts to give a tilting generator of $X$ or $X'$ over $Z$ as in \cite[Prop.~3.6]{ds12}, or it may be seen using the general machinery of \cite{dhl} and \cite{bfk}; compare especially \cite[Example~4.12]{dhl}.  Thus we get a \emph{window equivalence} $\psi_k := \iota'^* \circ (\iota^*)^{-1}$:
\[ \xymatrix{
& \W_k \ar[ld]_{\iota^*} \ar[rd]^{\iota'^*} \\
D^b(X) \ar[rr]_{\psi_k} && D^b(X').
} \]

\subsection*{Window equivalence vs.\ Bondal--Orlov equivalence}  The following was proved for $n=1$, $k=-1$, and $Z = \text{point}$ in \cite[Prop.~1]{ds13} 
by the same method, and it agrees with the analysis of general toric flops in \cite[\S3.1]{dhl_shipman}.

\begin{prop}
The equivalence $\psi_k$ coincides with $BO_{n+k}$.
\end{prop}
\begin{proof}
We will argue that the autoequivalence $\psi_k^{-1} \circ BO_{n+k}$ of $D^b(X)$ takes $\O_X(l)$ to itself for $k \le l \le k+n$, and acts as the identity on
\begin{equation} \label{sheaf_exts_to_act_as_1_on}
\E xt^i_{\varpi r}(\O_X(l), \O_X(l'))
\end{equation}
for $k \le l, l' \le k+n$ and $i \in \Z$.  Thus it acts as the identity on $\O_X(l) \otimes r^* \varpi^* \G$ for $k \le l \le k+n$ and $\G \in D^b(Z)$, and on all Exts between two such objects, hence on the category they split-generate, which is all of $D^b(X)$.

Clearly $\psi_k$ takes $\O_X(l)$ to $\O_{X'}(-l)$ for $k \le l \le k+n$.  Following the proof of \cite[Prop.~3.1]{kawamata2} we find that $BO_{n+k}$ does the same:
\begin{align*}
BO_{n+k}(\O_X(l))
&= p_*\!\left(\O_{\tilde X}((n+k)E) \otimes q^* \O_X(l)\right) \\
&= p_* \O_{\tilde X}(l-n-k,-n-k) \\
&= \O_{X'}(-l) \otimes p_* \O_{\tilde X}((n+k-l) E),
\end{align*}
and in the last term we have $p_* \O_{\tilde X}((n+k-l) E) = \O_{X'}$ because $0 \le n+k-l \le n$.  Thus $\psi_k^{-1} \circ BO_{n+k}$ takes $\O_X(l)$ to itself for $k \le l \le k+n$, and it remains to show that it acts as the identity on \eqref{sheaf_exts_to_act_as_1_on}.  This vanishes for $i \ne 0$ because $l'-l \ge -n$.  For $i=0$, we observe that on $X \setminus \P V$, both $\psi_k$ and $BO_{n+k}$ just act as the isomorphism $X \setminus \P V \cong X' \setminus \P V'$, so $\psi_k^{-1} \circ BO_{n+k}$ acts as the identity away from a set of codimension $n+1 \ge 2$, and because $\O_X(l)$ and $\O_X(l')$ are line bundles, a map between them is determined by its restriction to $X \setminus \P V$ by Hartogs' theorem.
\end{proof}

\subsection*{Window shift vs.\ spherical twist}

The following is a special case of \cite[Thm.~3.12]{ds12}, at least when $k=0$ and $Z = \text{point}$; some heuristic discussion of the case $n=1$ is given in [\emph{ibid}., \S2.2].  It also follows from the very general \cite[Prop.~3.4]{dhl_shipman}.

\begin{prop}
The autoequivalence $\psi_k^{-1} \circ \psi_{k+1}$ coincides with the spherical twist $T_k$ associated to the spherical functor
\[ F_k = j_*\big(\O_{\P V}(k) \otimes \varpi^*(-)\big)\colon D^b(Z) \to D^b(X) \]
of Theorem \ref{BO_family}.
\end{prop}
\begin{proof}
We sketch briefly the idea of the proof, taking $k=0$ for simplicity.  As in the previous proof, we let both functors act on generators $\O_X(1), \dotsc$, 
$\O_X(n), \O_X(n+1)$.  Clearly $\psi_0^{-1} \psi_1$ acts on $\O_X(1), \dotsc, \O_X(n)$ as the identity, and they are annihilated by the right adjoint
\[ R_0 = \varpi_* j^!(-) = \varpi_*\big(\O_{\P V}(-n-1) \otimes j^*(-)\big) \otimes \det V'[-n-1] \]
of $F_0$, so $T_0$ acts on them as the identity as well.  To understand $\O_X(n+1)$, we consider the Koszul resolution of the substack $[V \oplus 0\,/\,\C^*] \subset \XX$, which is cut out by a transverse section of $\O_\XX(-1) \otimes \varrho^* V'$:
\begin{multline*}
\O_\XX(n+1) \otimes \det V'^* 
\to \dotsb \\
\to \O_\XX(2) \otimes \Lambda^2 V'^*
\to \O_\XX(1) \otimes V'^*
\to \O_\XX \to \O_{[V \oplus 0\,/\,\C^*]}.
\end{multline*}
If we restrict to $X$, we get the Koszul resolution of $j(\P V) \subset X$:
\begin{multline} \label{koszul_of_zero_sec}
\O_X(n+1) \otimes \det V'^* 
\to \dotsb \\
\to \O_X(2) \otimes \Lambda^2 V'^*
\to \O_X(1) \otimes V'^*
\to \O_X \to j_* \O_{\P V}.
\end{multline}
If we restrict to $X'$, the last term goes away and we get $r'^*$ of the long Euler sequence of $\P V'$:
\begin{multline} \label{long_euler_seq}
\O_{X'}(-n-1) \otimes \det V'^* 
\to \dotsb \\
\to \O_{X'}(-2) \otimes \Lambda^2 V'^*
\to \O_{X'}(-1) \otimes V'^*
\to \O_{X'}.
\end{multline}
Now start with $\O_X(n+1)$ and apply $\psi_1$ to get $\O_{X'}(-n-1)$.  Use \eqref{long_euler_seq} to turn this into an $(n+1)$-term complex
\[ \underline{\O_{X'}(-n) \otimes V'} \to \O_{X'}(-n+1) \otimes \Lambda^2 V' \to \dotsb \to \O_{X'} \otimes \det V', \]
where the underlined term is in degree zero.  Apply $\psi_0^{-1}$ to get
\[ \underline{\O_X(n) \otimes V'} \to \O_X(n-1) \otimes \Lambda^2 V' \to \dotsb \to \O_X \otimes \det V'. \]
This is the middle $n+1$ terms of \eqref{koszul_of_zero_sec} tensored with $\det V'$, so rewrite it as an $(n+1)$-step extension of the last term by the first:
\[ \cone\left( j_* \O_{\P V} \otimes \det V'[-n-1] \to \O_X(n+1) \right). \]
On the other hand if we apply $T_0 = \cone(F_0 R_0 \to \id)$ to $\O_X(n+1)$ we get the same expression.  It remains to check that the two extensions are the same, and that $\psi_0^{-1} \psi_1$ and $T_0$ act in the same way on the Exts between $\O_X(1), \dotsc, \O_X(n+1)$, but for this we rely on the references given earlier.
\end{proof}

%% file: mukai.tex

\section{Mukai flops} \label{mukai}

Now assume the set-up of Theorem \ref{KN_family}: we have
\[ \xymatrixresized{
\P V \ar@{^{(}->}[r]^j \ar[d]_\varpi & X \\
{\phantom,}Z,
} \]
with $\N_{\P V/X} = \Omega^1_{\P V/Z}$ and $\HH^\text{odd}(Z) = 0$.

\begin{prop} \label{P-prop}
The functor
\[ F_k := j_*\big(\O_{\P V}(k) \otimes \varpi^*(-)\big)\colon D^b(Z) \to D^b(X) \]
is a $\P^n$-functor with $H = [-2]$.
\end{prop}
\begin{proof}
We emulate the proof of Theorem \ref{horja_thm} given on page \pageref{horja_pf_with_lb}.  The right adjoint of $F_k$ is
\[ R_k = \varpi_*\big(\O_{\P V}(-k) \otimes j^!(-)\big), \]
and we want to understand $R_k F_k$.  By Lemma \ref{j*j*}, the kernel inducing $j^! j_*$ has cohomology sheaves
\begin{equation} \label{wedge_of_normal}
\mathcal H^i = \Delta_* \Lambda^i \N_{\P V/X} = \Delta_* \Omega^i_{\P V/Z} \in \operatorname{Coh}(\P V \times \P V)
\end{equation}
and the monad structure (on the level of $\mathcal H^*$) is given by wedging.  To get the kernel inducing $\O_{\P V}(-k) \otimes j^! j_*(\O_{\P V}(k) \otimes -)$ we tensor with $\O_{\P V}(k) \boxtimes \O_{\P V}(-k)$, which does not change the sheaves \eqref{wedge_of_normal}.  To get the kernel inducing $R_k F_k$ we apply $(\varpi \times \varpi)_*$, which gives cohomology sheaves
\[ \mathcal H^i = \begin{cases}
\Delta_* \O_Z & i = 0, 2, \dotsc, 2n \\
0 & \text{otherwise},
\end{cases} \]
and the monad structure (again on the level of $\mathcal H^*$) is the right one.  But because $H\!H^\text{odd}(Z) = \Ext^\text{odd}_{Z \times Z}(\O_\Delta, \O_\Delta) = 0$, there are no extensions of $\O_\Delta[-2k]$ by $\O_\Delta[-2l]$, so the kernel inducing $R_k F_k$ splits as the sum of its cohomology sheaves. It is straightforward to verify that $R \cong H^n L$, and the proposition follows.
\end{proof}

\begin{cut}
\begin{rmk}
Alternatively we can deform to the normal bundle of $\P V$ in $X$.  On the central fiber the claim is proved for $k=0$ in \cite[\S3.2, Example~4]{nick}, and the proof for arbitrary $k$ is similar.  Now the object
\[ \O_\Delta \oplus \O_\Delta[-2] \oplus \dotsb \oplus \O_\Delta[-2n] \in D^b(Z \times Z) \]
is rigid because $H\!H^\text{odd}(Z) = 0$, and the monad structure condition of Definition \ref{P-def} is an open condition, so the claim holds on the general fiber.
\end{rmk}
\end{cut}

\begin{rmk} \label{rmk_on_HH_odd}
The assumption that $\HH^\text{odd}(Z) = 0$ is probably stronger than necessary, but it is satisfied in the applications we have in mind (where $Z$ is a K3 surface or a moduli space of sheaves on a K3 surface) and it greatly simplifies the proof of Proposition \ref{P-prop}.  \begin{cut}We could instead assume that the normal sequence of $\P V$ in $X$ splits, but this too is probably stronger than necessary, and it is hard to check.  Or we could adopt Cautis's definition of a $\P$-functor rather than Addington's, as discussed in Remark \ref{P-rmk}, and drop the requirement that $R_k F_k$ splits; but then we still have to require $\HH^1(Z) = 0$ in order to construct the $\P$-twist.\end{cut}
\end{rmk}

Now as in the introduction we let $X \xleftarrow{q} \tilde X \xrightarrow{p} X'$ be the Mukai flop of $X$ along $\P V$, we identify the exceptional divisor $E \subset \tilde X$ with the universal hyperplane in $\P V \times_Z \P V^*$, and we let
\[ \hat X = \tilde X \cup_E (\P V \times_Z \P V^*), \]
with maps $X \xleftarrow{\hat q} \hat X \xrightarrow{\hat p} X'$ given by $q$ and $p$ on $\tilde X$ and by the two projections on $\P V \times_Z \P V^*$.  Then we let $\L$ be the unique line bundle on $\hat X$ whose restriction to $\tilde X$ is $\O(E)$ and whose restriction to $\P V \times_Z \P V^*$ is $\O(-1,-1)$.
\begin{defn} \label{KN_k}
Let $X$, $X'$, and $\L \in \operatorname{Pic}(\hat X)$ be as in the previous paragraph.  For $k \in \Z$, we define
\[ \xymatrix{
D^b(X) \ar@<1ex>[rrrrrr]^*+{\KN_k := \hat p_*\big(\L^k \otimes \hat q^*(-)\big)} &&&&&&
D^b(X'). \ar@<1ex>[llllll]^*+{\KN'_k := \hat q_*\big(\L^k \otimes \hat p^*(-)\big)}
} \]
\end{defn}
\begin{thm}[{Kawamata \cite[\S5]{kawamata}, Namikawa \cite{namikawa}}]
The functors $\KN_k$ and $\KN'_k$ are equivalences.
\end{thm}
\noindent In fact Kawamata and Namikawa only prove this for $k=0$, but the generalization to arbitrary $k$ is straightforward.

To complete the proof of Theorem \ref{KN_family}, it remains to show that the $\P$-twist $P_k$ associated to the functor $F_k$ of Proposition \ref{P-prop} is isomorphic to $\KN'_{-k} \circ \KN_{n+k+1}$.  This will occupy the rest of the section.

\begin{prop} \label{BO_to_KN_in_local_model}
If $X$ is the total space of $\Omega^1_{\P V/Z}$ (the ``local model'') then $P_k = \KN'_{-k} \circ \KN_{n+k+1}$.
\end{prop}
\begin{proof}
The idea is that the Mukai flop is a hyperplane section of the standard flop, the Kawamata--Namikawa kernel is a hyperplane section of the Bondal--Orlov kernel, and the $\P$-twist is a hyperplane section of the spherical twist.

Let $\X$ be the total space of $\O_{\P V}(-1) \otimes \varpi^* V^*$.  Then the Euler sequence
\[ 0 \to \Omega^1_{\P V/Z} \to \O_{\P V}(-1) \otimes \varpi^* V^* \to \O_{\P V} \to 0 \]
determines a map $\X \to \A^1$ such that $X$ is the fiber over 0.  Following \cite[proof of Prop.~11.31]{huybrechts_fm}, we perform a standard flop $\X \leftarrow \tilde\X \rightarrow \X'$ along $\P V$ to get another family $\X' \to \A^1$ whose special fiber is $X'$.  Moreover the special fiber of $\tilde\X \to \A^1$ is $\hat X$, and the restriction of $\O_{\tilde\X}(kE)$ to $\hat X$ is $\L^k$, so we have
\begin{align} \label{BOi=iKN}
BO_k \circ \iota_* &= \iota'_* \circ \KN_k &
BO'_k \circ \iota'_* &= \iota_* \circ \KN'_k,
\end{align}
where $\iota\colon X \to \X$ and $\iota'\colon X' \to \X'$ are the inclusions.

Now $\iota_* \circ F_k$ is the spherical functor studied in \S\ref{sod}\begin{cut} and \S\ref{windows}\end{cut}; writing $T_k$ for the associated spherical twist, we have $T_k = BO'_{-k} \circ BO_{n+k+1}$ by Theorem \ref{BO_family}.  We claim next an isomorphism of functors
\begin{equation} \label{TiP=i}
T_k \circ \iota_* \circ P_k^{-1} = \iota_*.
\end{equation}
By Proposition \ref{intertwinement_prop}, their kernels in $D^b(X \times \X)$ agree on the fiber over each point of $Z$. In particular the left-hand functor takes skyscrapers sheaves of points $\O_x$ to skyscraper sheaves $\O_{\iota(x)}$, so by \cite[Cor.~1.12]{bbh} its kernel is a line bundle $\M$ on the graph of $\iota$.  This $\M$ is trivial on the fiber over each point of $Z$, hence is pulled back from $Z$.  But $T_k$ is the identity away from $\P V \subset \X$, and $P_k$ is the identity away from $\P V \subset X$, so $\M$ is furthermore trivial away from $\P V \subset X$, which has codimension $n$.  If $n \ge 2$ we conclude that $\M$ is trivial.  If $n=1$ we conclude that $\M \cong \O_X(m \P V)$ for some $m \in \Z$; but the latter is not pulled back from $Z$ unless $m=0$, because $\O_X(\P V)|_{\P V} = \N_{\P V/X} = \omega_{\P V/Z} = \O_{\P V}(-2) \otimes \varpi^* \det V^*$, so $\M$ is again trivial.  Thus we have established \eqref{TiP=i}.

Now combining \eqref{BOi=iKN} and \eqref{TiP=i} we get
\begin{align}
\iota_* \circ \KN_{-k}' \circ \KN_{n+k+1} \circ P_k^{-1} = \iota_*.
\end{align}
To finish the proof, it is enough to observe that any Fourier--Mukai functor $\Phi\colon D^b(X) \to D^b(X)$ is isomorphic to $\id_X$ if and only if $\iota_* \circ \Phi = \iota_*$.  Thinking about kernels, this says that an object $\Phi \in D^b(X\times X)$ is isomorphic to $\O_\Delta$ if and only if $(1 \times \iota)_* \Phi$ is isomorphic to $(1 \times \iota)_* \O_\Delta$.  Since $(1 \times \iota)_*$ is exact we see that $(1 \times \iota)_* \Phi$ is a sheaf (rather than a complex) if and only if $\Phi$ is; and pushing forward by a closed embedding is fully faithful on sheaves (though of course not on complexes) so $(1 \times \iota)_* \Phi \cong (1 \times \iota)_* \O_\Delta$ if and only if $\Phi \cong \O_\Delta$.
\end{proof}

\begin{cut}
\begin{rmk}
The same proof applies when $X$ is compact hyperk\"ahler, which is the main source of examples: following \cite[Rmk.~11.32]{huybrechts_fm} we observe that by \cite[Lem.~3.6]{huybrechts_non_sep} there is a family $\X$ over a curve such that $\X_0 = X$ and $\N_{\P V/\X} = \O_{\P V}(-1) \otimes \varpi^* V^*$.
\end{rmk}
\end{cut}

\begin{prop} \label{defo_to_normal_cone}
If $X$ is arbitrary then $P_k = \KN'_{-k} \circ \KN_{n+k+1}$.
\end{prop}
\begin{proof}
We will show that $\Phi := \KN_{-k}' \circ \KN_{n+k+1} \circ P_k^{-1}$ is the identity.  First we will argue that $\Phi(\O_x) = \O_x$ for all points $x \in X$; this is immediate when $x \in X \setminus \P V$, so the interesting case is when $x \in \P V$.  Thus $\Phi$ is given by tensoring by a line bundle, which moreover is trivial away from $\P V \subset X$, and this has codimension $n$; if $n \ge 2$ we conclude that the line bundle is trivial, and if $n=1$ we conclude that it is $\O_X(m \P V)$ for some $m \in \Z$.  But the restriction of the latter to a fiber of $\varpi\colon \P V \to Z$ is $\O_{\P^1}(-2m)$, so we will argue that $\Phi$ takes $\O_{j(\varpi^{-1}(z))} = \O_{\P^1}$ to itself for $z \in Z$, hence $m=0$.

First then we fix $x \in \P V$ and argue that $\Phi(\O_x) = \O_x$.  Consider the deformation to the normal bundle of $\P V$ in $X$ \cite[Ch.~5]{fulton}: thus we have a family $\X \to \A^1$ such that $\X_0$ is the total space of $\N_{\P V/X} = \Omega^1_{\P V/Z}$, and $\X_t = X$ for all $t \ne 0$.  The sheaves inducing $F_k$ and $\KN_k$ can be constructed flatly in the family, so we get an object $\G \in D^b(\X)$ such that $\G_t := \G \otimes \O_{\X_t}$ is the appropriate $(\KN_{-k}' \circ \KN_{n+k+1} \circ P_k^{-1})(\O_x)$ for $\X_t$: thus $\G_0 = \O_x$ by the previous proposition, and $\G_t = \Phi(\O_x)$ for $t \ne 0$.  We will argue that $\G$ is a sheaf (rather than a complex) supported on $x \times \A^1 \subset \X$ and flat over $\A^1$, so $\G_t = \O_x$ for all $t$.

Let $S_q = \supp(\H^q(\Phi(\O_x)) \subset \P V \subset X$; we will find that $S_q = \emptyset$ for $q \ne 0$ and $S_0 = \{x\}$.  Because $\G_t = \Phi(\O_x)$ for $t \ne 0$, we see that $\supp(\H^q(\G))$ contains
\[ S_q \times (\A^1 \setminus 0) \ \subset\ \P V \times \A^1 \ \subset\ \X. \]
But $\supp(\H^q(\G))$ is closed, so in fact it contains $S_q \times \A^1$.  Now use the Grothendieck spectral sequence
\begin{equation} \label{groth_spec_seq}
E_2^{p,q} = \Tor_{-p}(\H^q(\G), \O_{\X_0}) \Rightarrow \H^{p+q}(\G_0).
\end{equation}
Because $\X_0$ is a divisor in $\X$, these Tors vanish apart from $\Tor_0$ and $\Tor_1$, and the sequence degenerates at the $E_2$ page\begin{cut}, displayed in Figure \ref{E2_page} below\end{cut}.

\begin{cut}
\begin{figure}
\[ \xymatrix@R=3pt@C=20pt{
 & \vdots & \vdots \\
0 \ar[rrdd] & \Tor_1(\H^1(\G), \O_{\X_0}) \ar[rrdd] & \Tor_0(\H^1(\G), \O_{\X_0}) & 0 \\ \\
0 \ar[rrdd] & \Tor_1(\H^0(\G), \O_{\X_0}) \ar[rrdd] & \Tor_0(\H^0(\G), \O_{\X_0}) & 0 \\ \\
0 & \Tor_1(\H^{-1}(\G), \O_{\X_0}) & \Tor_0(\H^{-1}(\G), \O_{\X_0}) & 0 \\
 & \vdots & \vdots
}  \]
\caption{$E_2$ page of the spectral sequence \eqref{groth_spec_seq}.} \label{E2_page}
\end{figure}
\end{cut}

Since $\G_0 = \O_x[0]$, we see that for $q \ne 0$ we have $\Tor_0(\H^q(\G), \O_{\X_0}) = 0$, so $\supp(\H^q(\G)) \cap \X_0 = \emptyset$, so $S_q = \emptyset$, so $\H^q(\G) = 0$; thus $\G$ is a sheaf.  Next we see that $\Tor_1(\H^0(\G), \O_{\X_0}) = 0$, so that sheaf is flat over $\A^1$.  Finally we see that $\Tor_0(\H^0(\G), \O_{\X_0}) = \O_x$, so $\supp(\H^0(\G)) \cap \X_0 = \{x\}$, so $S_0 = \{x\}$, and we conclude that $\G_t = \O_x$ for all $t$, so $\Phi(\O_x) = \O_x$ as desired.

The argument that $\Phi$ takes $\O_{j(\varpi^{-1}(z))} = \O_{\P^n}$ to itself is entirely similar: we get a sheaf $\G'$ on $\X$ supported on $\P^n \times \A^1$, flat over $\A^1$, and having $\G'_0 = \O_{\P^n}$; the latter is rigid, so $\Phi(\O_{\P^n}) = \G'_{t\ne 0} = \O_{\P^n}$.
\end{proof}

\begin{cut}
\begin{rmk}
We could instead have deduced Proposition \ref{defo_to_normal_cone} from Proposition \ref{BO_to_KN_in_local_model} using Kawamata and Namikawa's observation that if $Z$ is a point then the formal neighborhood of $\P V$ in $X$ is isomorphic to the formal neighborhood of the zero section in the total space of $\Omega^1_{\P V}$;
%
%
%
%
if $Z$ is general then this may not hold,
%
%
but we can restrict to an open set in $Z$ over which $V$ is trivial and go from there as in \cite[\S5]{namikawa}.  We have used deformation to the normal bundle partly for the sake of novelty, but also because it is a simpler and more widely-applicable technique.
\end{rmk}
\end{cut}